\documentclass[12pt]{article}

\usepackage{amsmath}
\usepackage{amsthm}
\usepackage{amssymb}
\usepackage{graphicx}
\usepackage{color}
\usepackage{bm}

\allowdisplaybreaks[4]

\numberwithin{equation}{section}
\newtheorem{lm}{Lemma}
\newtheorem{thm}{Theorem}
\newtheorem{exam}{Example}
\title{\bf
Continuous solutions of a second order iterative equation
\thanks{Supported by NSFC grants \# 11521061 and \# 11231001 and by PCSIRT IRT\_15R53
}
}

\author{
{\normalsize
{\sc Xiao Tang}~~~and ~~~
{\sc Weinian Zhang}
\thanks{Corresponding author: matzwn@126.com
}
       }
\\
{\small  Department of Mathematics, Sichuan University}
\\
{\small  Chengdu, Sichuan 610064, China}
}

\date{}

\begin{document}
	
\maketitle
	
\begin{abstract}
In this paper we study the existence of continuous solutions and their constructions 
for a second order iterative functional equation
which involves iterate of the unknown function 
and a nonlinear term.
Imposing Lipschitz conditions to those given functions, we prove 
the existence of continuous solutions on the whole $\mathbb{R}$ by applying the contraction principle.
In the case without Lipschitz conditions we hardly use the contraction principle,
but we construct continuous solutions on $\mathbb{R}$ recursively with a partition of $\mathbb{R}$.

\indent
\vskip 0.1cm
		
		
\vskip 0.1cm
		
		
\end{abstract}



\section{Introduction}

In the problem session of the 38th ISFE held in 2000 in Hungary, 
 N. Brillou\"{e}t-Belluot (\cite{Belluot}) proposed the second order iterative equation
\begin{eqnarray}\label{speciallinear}
\varphi^2(x)=\varphi(x+a)-x,~~~~  x\in \mathbb{R},
\end{eqnarray}
and asked: {\it What are its solutions?}
Three years later K. Baron (\cite{Baron}) emphasized it again.
This equation was reduced from the multi-variable equation
\begin{eqnarray*}
x+\varphi(y+\varphi(x))=y+\varphi(x+\varphi(y)),
\end{eqnarray*}
an important functional equation 
which has been attractive to many researchers
(\cite{Sablik, Yarczyk, Balcerowski}).
For the case $a=0$, equation \eqref{speciallinear} has no continuous solutions by  Theorem $11$ of (\cite{Matkowski}) or Theorem $5$ of (\cite{YDL}).
For the case $a \ne 0$, by Theorem $1$ of \cite{ZYY} equation \eqref{speciallinear} also has no continuous solutions.

In $2010$ N. Brillou\"{e}t-Belluot and W. Zhang (\cite{ZWN}) 
investigated a more general form 
\begin{eqnarray}\label{linear}
\varphi^2(x)=\lambda \varphi(x+a) +\mu x,
\end{eqnarray}
where $\lambda, a$ and $\mu$ are all real such that $a\lambda \ne 0$. They
used the contraction principle 
to prove the existence of a continuous solution
under the condition
\begin{eqnarray}\label{ZWN}
|\lambda| > \max \{2,2\sqrt{2|\mu|}\}   ~~\text{or}~~ 1+2|\mu| < |\lambda| \le 2
\end{eqnarray}
and employed the technique of piecewise construction to obtain piecewise continuous
solutions in the case that
$0\le \mu<1$ and $\lambda\ge 2(1-\mu)$.
Later Y. Zeng and W. Zhang (\cite{ZYY}) proved that equation \eqref{linear} has no continuous
solutions on $\mathbb{R}$
if $\lambda=1$ and $\mu \le -1$, 
which is the source result that implies the nonexistence 
stated in the end of last paragraph.
They also gave existence of continuous solutions on $\mathbb{R}$
in the case that
\begin{eqnarray}\label{ZYY1}
|\lambda| \in (2,+\infty)~~{\rm and}~~  \mu \in [-\lambda^2/4,\lambda^2/4]
\end{eqnarray}
and the case that
\begin{eqnarray}\label{ZYY2}
|\lambda| \in (1,2] ~~{\rm and}~~ \mu \in (1-|\lambda|,|\lambda|-1).
\end{eqnarray} 

In this paper we generally consider the iterative equation
\begin{eqnarray}\label{rose}
\varphi^2(x)= h(\varphi(f(x)))+g(x), ~~~~ x\in \mathbb{R},
\end{eqnarray}
where $h, f$ and $g: \mathbb{R}\to \mathbb{R}$ 
are given continuous functions and $\varphi:\mathbb{R}\to \mathbb{R}$ is the unknown one.
This equation includes equation $\eqref{linear}$ as a special case with 
the choice that $f(x)=x+a,~h(x)=\lambda x$ and $g(x)=\mu x$. 
In section 2 we consider bounded $g$ and prove the existence of a bounded continuous solution on $\mathbb{R}$ (Theorem \ref{thm1}) 
under Lipschitz conditions to those given functions or their inverses
by applying the contraction principle.
Section 3 is devoted to the case of unbounded $g$. 
We give a result of the existence (Theorem \ref{thm2}) on compact intervals by 
modifying Theorem \ref{thm1} 
and obtain another result of the existence (Theorem \ref{thm3}) on the whole $\mathbb{R}$ 
with additional assumptions of bounded nonlinearities by applying the contraction principle.
In section 4 we discuss equation (\ref{rose}) in the case without Lipschitz conditions, where
we hardly apply the contraction principle again.
We construct continuous solutions recursively with a partition of $\mathbb{R}$
in some cases (Theorem \ref{thm4}).
We finish this paper in section 5 with some remarks.


\section{Case of bounded $g$}

We need the following hypotheses:
\begin{description}
		\item {(\bf C1)} $h$ is uniformly expansive, i.e., there exists a constant $K > 1$ such that
		\begin{eqnarray}\label{expansive}
			| h(x)-h(y) | \ge K | x-y |,~~~~ \forall x,y \in \mathbb{R};
		\end{eqnarray}

\item {(\bf C2)} There is a constant $\alpha >0$ such that
			\begin{eqnarray}\label{ILip}
			| f(x) -f(y) | \ge \alpha | x-y |, ~~~~ \forall x,y\in \mathbb{R};
			\end{eqnarray}

\item {(\bf C3)} $g \in C^0_b(\mathbb{R})$ is Lipschitzian with the Lipschitz constant $\text{Lip}(g)\le \beta$.
\end{description}

	
\begin{thm}
Suppose that functions $h,f$ and $g$ fulfill conditions ${\bf (C1)}$, ${\bf (C2)}$ and ${\bf (C3)}$, where constants $K, \alpha$ and $\beta$ 
satisfy 
\begin{eqnarray}
&&\beta \le \frac{1}{4}\alpha^2K^2 ~~~~~~~~~~~~~~~~\,~~~~ \text{when}~~ \alpha  < 2(1-\frac{1}{K}),
\label{contracself1}
\\
&&\beta < (K-1)(\alpha K-K+1)~~ \text{when}~~ \alpha  \ge 2(1-\frac{1}{K}). 
\label{contracself2}
\end{eqnarray}
Then equation $\eqref{rose}$ has a bounded continuous solution on $\mathbb{R}$.
\label{thm1}
\end{thm}


\noindent
{\bf Proof.
} 
From assumption ${\bf (C2)}$ we get $f: \mathbb{R}\to \mathbb{R}$ is a homeomorphism. In fact,
it is clear that $f$ is injective. In order to prove ``onto'' for $f$, without loss of generality we assume that 
$A:=\lim_{x\to +\infty}f(x)$ exists. Choosing $y=0$ and letting $x\to+\infty$ in $\eqref{ILip}$,
we get a contradiction.
Moreover, 
assumption ${\bf (C1)}$ implies that $h$ is bijective and its inverse $h^{-1}$ is contractive.
Actually, the method to prove bijection of $h$ is same as that for $f$. 
Since $K>1$, inequality $\eqref{expansive}$ yields that $h^{-1}$ is contractive. 
Thus, under conditions ${\bf (C1)}$ and ${\bf (C2)}$ equation $\eqref{rose}$ is equivalent to the form
\begin{eqnarray}\label{peony}
	\varphi(x)=h^{-1}(\varphi^2 \circ f^{-1}(x)-g\circ f^{-1}(x)),~~~~ x\in \mathbb{R},
\end{eqnarray}
where $\circ$ denotes the composition of functions, i.e., $f\circ g(x):=f(g(x))$.

Clearly, the set   
$$
C^0_b(\mathbb{R}):=\{\varphi: \mathbb{R} \to \mathbb{R} | \varphi~ \text{is continuous and}~
	 \underset{x\in \mathbb{R}}\sup \lvert \varphi(x) \lvert < +\infty \}
$$
is a Banach space equipped with the norm 
$
\|\varphi\|:=
\sup_{x\in \mathbb{R}}|\varphi(x)|.
$ 	 
Let 
$
C^0_b(\mathbb{R};L):= C^0_b(\mathbb{R}) \cap \{ \varphi:\mathbb{R} \to \mathbb{R} | 
\text{Lip} (\varphi) \leq L \}, 
$
where $L>0$ is a constant.
Clearly, $C^0_b(\mathbb{R};L)$ is a closed subset of $C^0_b(\mathbb{R})$.
Define a mapping $\mathcal{T}: C^0_b(\mathbb{R};L) \to C^0_b(\mathbb{R})$ such that for a given function $\varphi \in C^0_b(\mathbb{R};L)$
\begin{eqnarray}
\mathcal{T}\varphi(x)= h^{-1}(\varphi^2\circ f^{-1}(x)-g \circ f^{-1}(x)).
\label{TTT}
\end{eqnarray}
Clearly, functions $\varphi$ is a solution of equation $\eqref{peony}$ if and only if 
$\varphi$ is a fixed point of the mapping $\mathcal{T}$. 
For a given function $\varphi \in C^0_b(\mathbb{R};L)$,
it is obvious that $\mathcal{T}\varphi$ is a continuous function. 
Since $\varphi$ is bounded,
let $M_*:=\max \{\|\varphi \|,\|g\|\}$. Then $\|\varphi^2-g\|\le \|\varphi\|+\|g\|\le 2M_*$.
It follows that
\begin{eqnarray*}
\begin{split}
	\underset{x\in \mathbb{R}}{\sup}|\mathcal{T}\varphi(x) |    
	&  =    \underset{x\in \mathbb{R}}{\sup} |h^{-1}(\varphi^2\circ f^{-1}(x)-g\circ f^{-1}(x))|
	\\
	&  =    \underset{x\in \mathbb{R}}{\sup} |h^{-1}(\varphi^2(x)-g(x))|
	\\
	&  \le    \underset{|x|\le 2M_*}{\sup} |h^{-1}(x)|<+\infty,
\end{split}	
\end{eqnarray*}
implying that $\mathcal{T}\varphi\in C_b^0(\mathbb{R})$, i.e., $\mathcal{T}$ given in (\ref{TTT}) is well defined.

We claim that $\mathcal{T}$ is a self-mapping on $C^0_b(\mathbb{R};L)$
for an appropriate constant $L>0$. 
In fact, for any $x_1,x_2\in \mathbb{R}$,
\begin{eqnarray*}
&&|\mathcal{T}\varphi(x_1)-\mathcal{T}\varphi(x_2)| 
\\
&&\quad = 
|h^{-1}(\varphi^2 \circ f^{-1}(x_1)-g \circ f^{-1}(x_1))- h^{-1}(\varphi^2 \circ f^{-1}(x_2)-g \circ f^{-1}(x_2))|
\\
&&\quad \leq 
\frac{1}{K} |\varphi^2\circ f^{-1}(x_1)-g\circ f^{-1}(x_1)-\varphi^2\circ f^{-1}(x_2)+g\circ f^{-1}(x_2)|
\\
&&\quad 
\leq \frac{1}{K} (\frac{L^2}{\alpha}+\frac{\beta}{\alpha})|x_1-x_2|
\\
&&\quad 
\leq L|x_1-x_2|
\end{eqnarray*}
if $L$ is chosen such that
$$
\frac{1}{K} (\frac{L^2}{\alpha}+\frac{\beta}{\alpha})\le L,
$$
which is equivalent to the inequalities
\begin{eqnarray}\label{self}
\frac{1}{2}\alpha K-\frac{1}{2}\sqrt{\alpha^2 K^2-4\beta} \le L \le \frac{1}{2}\alpha K+\frac{1}{2}\sqrt{\alpha^2 K^2-4\beta},
\end{eqnarray}
and
\begin{eqnarray}\label{Delta}
\beta\le \frac{\alpha^2 K^2}{4}.
\end{eqnarray}

Next, we assert that $\mathcal{T}$ is a contraction on $C^0_b(\mathbb{R};L)$
if constant $L$ satisfies
\begin{eqnarray}\label{contrac}
L < K-1
\end{eqnarray}
because for any two functions $\varphi_1, \varphi_2 \in C^0_b(\mathbb{R};L)$,
\begin{eqnarray*}
	\begin{split}
		\| \mathcal{T}\varphi_1 -\mathcal{T}\varphi_2 \| &= \underset{x\in\mathbb{R}}{\sup}|\mathcal{T}\varphi_1(x) -\mathcal{T}\varphi_2(x)|
		\\
		&=\underset{x\in\mathbb{R}}{\sup}|h^{-1}(\varphi_1^2\circ f^{-1}(x)-g\circ f^{-1}(x))-
		\\
		&~~~~	h^{-1}(\varphi_2^2\circ f^{-1}(x)-g\circ f^{-1}(x))|
		\\
		&\leq \frac{1}{K} \underset{x\in \mathbb{R}}{\sup}|\varphi_1^2\circ f^{-1}(x)-\varphi_2^2\circ f^{-1}(x)|
		\\
		&= \frac{1}{K} \underset{x\in\mathbb{R}}{\sup}|\varphi_1^2(x)-\varphi_2^2(x)|
		\\
		&\leq \frac{1}{K}\{\underset{x\in \mathbb{R}}{\sup}|\varphi_1^2(x)-\varphi_1(\varphi_2(x))|+
		\underset{x\in \mathbb{R}}{\sup}|\varphi_1(\varphi_2(x))-\varphi^2_2(x)|\}
		\\
		&\leq \frac{1}{K}(L+1)\|\varphi_1-\varphi_2\|.
	\end{split}
\end{eqnarray*}

As a consequence, the above claim and assertion conclude that
$\mathcal{T}$ is a contractive self-mapping on $C^0_b(\mathbb{R};L)$ 
if we can choose $L$ to fulfill $\eqref{self}$ and $\eqref{contrac}$ under  (\ref{Delta}).
Note that 
such $L$ exists if
\begin{eqnarray}\label{intersection}
\frac{1}{2}\alpha K-\frac{1}{2}\sqrt{\alpha^2 K^2-4\beta} < K-1.
\end{eqnarray}
Clearly,
when $\alpha < 2(1-\frac{1}{K})$, inequality $\eqref{intersection}$ holds automatically and 
thus we only need condition (\ref{Delta}), which is the same as   
condition $\eqref{contracself1}$.
When $\alpha \ge 2(1-\frac{1}{K})$, 
inequality $\eqref{intersection}$
is simplified to $\beta < (K-1)(\alpha K-K+1)$. 
Associated with (\ref{Delta}), it requires 
\begin{eqnarray}\label{draft2}
\beta < \min\{(K-1)(\alpha K-K+1), \frac{1}{4} \alpha^2 K^2\} ~~\text{as}~~ \alpha \ge 2(1-\frac{1}{K}).
\end{eqnarray}
Since 
$$
(K-1)(\alpha K-K+1)\le(\frac{K-1+\alpha K -K+1}{2})^2 = \frac{1}{4} \alpha^2 K^2,
$$
$\eqref{draft2}$
becomes
$\beta < (K-1)(\alpha K-K+1)$ when $\alpha \ge 2(1-\frac{1}{K})$, which is the same as condition $\eqref{contracself2}$.
Summarily,
under condition $\eqref{contracself1}$ or $\eqref{contracself2}$ we can choose an appropriate constant $L$ such that
$
\frac{1}{2}\alpha K -\frac{1}{2} \sqrt{\alpha^2 K^2-4\beta} \le L < K-1,
$
which guarantees $\eqref{self}$ and $\eqref{contrac}$ to be true.
By the Contraction Principle,
$\mathcal{T}$ 
has a unique fixed point in $C^0_b(\mathbb{R};L)$, which gives a solution.
The proof is completed.
\hfill$\square$

\begin{figure}
	\begin{minipage}[t]{0.5\linewidth}
		\centering
		\includegraphics[width=2.2in]{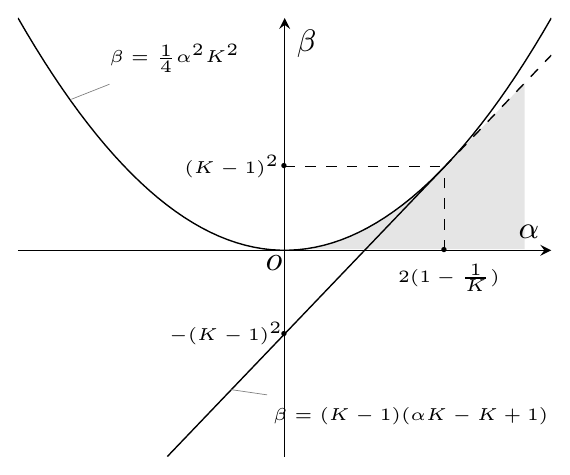}
		\caption{Region of $(\alpha, \beta)$ for \eqref{contracself1} and \eqref{contracself2}.}
		\label{shadow}
	\end{minipage}
	\begin{minipage}[t]{0.5\linewidth}
		\centering
		\includegraphics[width=2.55in]{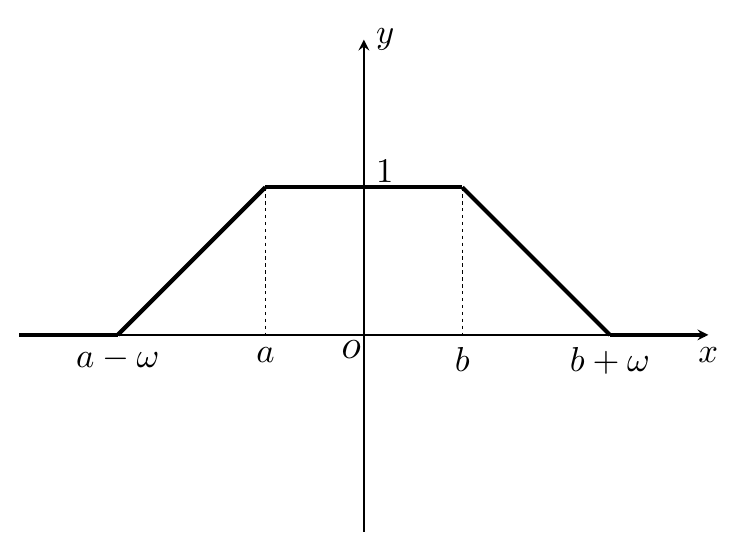}
		\caption{Graph of $\sigma_{\omega}$.}
		\label{omega}
	\end{minipage}
\end{figure}

As shown in Figure \ref{shadow},
for each given $K>1$
conditions (\ref{contracself1}) and (\ref{contracself2}) hold in the left shadowed region 
and the right shadowed region of the dashed line $x=2(1-\frac{1}{K})$ in the $(\alpha, \beta)$-plane respectively, 
from which we easily choose two examples: one is that $K=2$, $\alpha=\frac{1}{2}$, $\beta=\frac{1}{8}$ and 
the other is that
$K=\alpha=\beta=2$, which satisfy $\eqref{contracself1}$ and $\eqref{contracself2}$ respectively.

\begin{exam}{\rm
Our Theorem \ref{thm1} is applicable to the equation
\begin{eqnarray}\label{exam1}
\varphi^2(x)=2\varphi(2x) +\sin x,~~~~  x\in \mathbb{R},
\end{eqnarray}
which is of the form (\ref{rose}), where $f(x)=h(x)=2x$ and $g(x)=\sin x$. One can check that
$f, g$ and $h$ satisfy conditions {\bf (C1)}-{\bf (C3)} with constants
$K=\alpha=\beta=2$. Further, 
$2(1-1/K)=1<\alpha$. Thus, we can verify that 
$$
(K-1)(\alpha K -K +1)=3>\beta,
$$ 
i.e., 
condition (\ref{contracself2}) is fulfilled. By our Theorem \ref{thm1},
equation $\eqref{exam1}$ has a bounded continuous solution on $\mathbb{R}$. 
}
\end{exam}


\section{Case of unbounded $g$}

Theorem \ref{thm1} requires $g\in C^0_b(\mathbb{R})$, i.e., $g$ is a bounded function.
With a modification, we can obtain the following Theorem for unbounded $g$ but
the solution is not defined on the whole $\mathbb{R}$.


\begin{thm}
Suppose that conditions {\bf (C1)} and {\bf (C2)} hold and $g$ is Lipschitzian on $\mathbb{R}$
with Lip$(g)\le \beta$,
where $\beta$ satisfies \eqref{contracself2} when $\alpha \ge 2(1-\frac{1}{K})$ and
	\begin{eqnarray}
	&&\beta < \frac{1}{4}\alpha^2K^2 ~~~~~~~~ \text{when}~~\alpha < 2(1-\frac{1}{K}).
	\label{contracselfU1}
	\end{eqnarray}
Then for any compact interval $I$ equation $\eqref{rose}$ has a continuous solution on $I$.
\label{thm2}
\end{thm}


Before proving Theorem \ref{thm2}, we make a truncation to the function $g$. 
For a given compact interval $I=[a,b]$ and a number $\omega>0$, 
consider the function
\begin{eqnarray}\label{ccc}
\sigma_{\omega}(x):=
\begin{cases}
1,~~~~  & x\in I,
\\
\frac{1}{\omega}x+1-\frac{a}{\omega},~~~~  & x\in (a-\omega,a),
\\
-\frac{1}{\omega}x+1+\frac{b}{\omega},~~~~ & x\in (b,b+\omega),
\\
0,~~~~      & x\in (-\infty,a-\omega] \cup [b+\omega,+\infty),
\end{cases}
\end{eqnarray}
as shown in Figure $2$,
which is Lipschitzian with ${\rm Lip}(\sigma_\omega) = 1/\omega$.
Let
\begin{eqnarray}\label{cutoff}
\tilde{g}(x) := g(\sigma_{\omega}(x)x),~~~~ x\in \mathbb{R}.
\end{eqnarray}
One can check that $\tilde{g}$ is bounded and  
\begin{eqnarray}\label{cut}
\tilde{g}(x)=g(x),~~\forall x\in I.
\end{eqnarray}
For the estimation of Lipschitzian constants of $\tilde{g}$, we have the following lemma.


\begin{lm}
	Suppose that functions $g$ is Lipschitzian on $\mathbb{R}$ with  Lip$(g)\le \beta$. Then 
	\begin{eqnarray}
	\text{Lip}(\tilde{g}) \le \beta(1+\frac{\max\{|a|, |b|\}}{\omega}),
	\end{eqnarray}
	where the function $\tilde{g}$ is defined as \eqref{cutoff}.
	\label{lm1}
\end{lm}

\noindent
{\bf Proof.}
We need prove that for arbitrary $x_1,x_2 \in \mathbb{R}$
\begin{eqnarray}
|\tilde{g}(x_1)-\tilde{g}(x_2)| \le \beta(1+\frac{\max\{|a|, |b|\}}{\omega})|x_1-x_2|.
\end{eqnarray}
First, consider the case that one of $x_1,x_2$ does not belong to $(a-\omega, b+\omega)$, without loss of generality, assume that $x_1 \notin (a-\omega, b+\omega)$. According to the definition of the function $\sigma_{\omega}$, it follows that $\sigma_{\omega}(x_1) =0$.
Then
\begin{eqnarray*}
	|\tilde{g}(x_1)-\tilde{g}(x_2)|&=&|g(0)-g(\sigma_{\omega}(x_2)x_2)|
	\\
	&\le& \beta|\sigma_{\omega}(x_2)x_2|
	\\
	&=&
	\beta|\sigma_{\omega}(x_2)x_2-\sigma_{\omega}(x_1)x_2| 
	\\
	&\le& \beta(\frac{\max\{|a-\omega|, |b+\omega|\}}{\omega})|x_1-x_2|
	\\
	& \le& \beta(1+\frac{\max\{|a|,|b|\}}{\omega})|x_1-x_2|. 
\end{eqnarray*}
Next, we divide the opposite case that $x_1,x_2 \in (a-\omega, b+\omega)$ into two subcases. Case {\bf (A)}: $x_1x_2>0$. Case {\bf (B)}: $x_1x_2\le 0$.
For case {\bf (A)}, without loss of generality, assume that
$0 < x_1 < x_2$.
We discuss case {\bf (A)} by the following subcases. Case {\bf (A1)}: one of $x_1,x_2$ belongs to $(a-\omega, a)$; case {\bf (A2)}: $x_1,x_2 \in [a,b+\omega)$. For case {\bf (A1)} we assume that 
$x_1 \in (a-\omega, a)$
without loss of generality. Then
\begin{eqnarray*}
	|\tilde{g}(x_1)-\tilde{g}(x_2)|&=&|g(\sigma_{\omega}(x_1)x_1)-g(\sigma_{\omega}(x_2)x_2)|
	\\
	&\le& \beta|\sigma_{\omega}(x_1)x_1- \sigma_{\omega}(x_2)x_2|
	\\
	&=&
	\beta(|\sigma_{\omega}(x_1)x_1- \sigma_{\omega}(x_2)x_1| + |\sigma_{\omega}(x_2)x_1- \sigma_{\omega}(x_2)x_2| )
	\\
	&\le& \beta( \frac{x_1}{\omega} +1)|x_1-x_2|
	\\
	& \le& \beta(\frac{a}{\omega} +1)|x_1-x_2|
	\\
	& \le& \beta(1+\frac{\max\{|a|,|b|\}}{\omega})|x_1-x_2| 
\end{eqnarray*}
because of $0< x_1 <a$. For case {\bf (A2)}, noting that  $0< x_1<x_2$,
$\sigma_{\omega}(x_1) >0$ and $0< \sigma_{\omega}(x_2) \le \sigma_{\omega}(x_1)$,
we have
$$
0<\sigma_{\omega}(x_1)x_1 < \sigma_{\omega}(x_1)x_2 ~~~ \text{and} ~~~0< \sigma_{\omega}(x_2)x_2 \le \sigma_{\omega}(x_1)x_2,
$$
which imply that 
\begin{eqnarray*}
	|\tilde{g}(x_1)-\tilde{g}(x_2)|&=&|g(\sigma_{\omega}(x_1)x_1)-g(\sigma_{\omega}(x_2)x_2)|
	\\
	&\le& \beta|\sigma_{\omega}(x_1)x_1- \sigma_{\omega}(x_2)x_2|
	\\
	&\le&\beta | \sigma_{\omega}(x_1)x_1 - \sigma_{\omega}(x_1)x_2 | (\text{or} ~\beta | \sigma_{\omega}(x_2)x_2 - \sigma_{\omega}(x_1)x_2 |)
	\\
	&\le& \beta |\sigma_{\omega}(x_1)||x_1-x_2| ( \text{or} ~ \frac{\beta x_2}{\omega} |x_1 -x_2 | )
	\\
	&\le & \beta |x_1-x_2| (\text{or} ~ \frac{\beta (b+\omega)}{\omega} | x_1-x_2 |)
	\\
	& \le& \beta(1+\frac{\max\{|a|,|b|\}}{\omega})|x_1-x_2|. 
\end{eqnarray*}
In case {\bf (B)}, without loss of generality, we assume that $x_1\le 0 < x_2$. 
Noticing that $0< \sigma_{\omega}(x_1), \sigma_{\omega}(x_2) < 1$ as $x_1, x_2 \in (a- \omega, b+\omega)$, we acquire that
$|\sigma_{\omega}(x_1)x_1- \sigma_{\omega}(x_2)x_2| \le |x_1-x_2|$, which implies that
\begin{eqnarray*}
	|\tilde{g}(x_1)-\tilde{g}(x_2)|&=&|g(\sigma_{\omega}(x_1)x_1)-g(\sigma_{\omega}(x_2)x_2)|
	\\
	&\le& \beta|\sigma_{\omega}(x_1)x_1- \sigma_{\omega}(x_2)x_2|
	\\
	&\le & \beta |x_1-x_2|
	\\
	& \le& \beta(1+\frac{\max\{|a|,|b|\}}{\omega})|x_1-x_2|. 
\end{eqnarray*}
Thus, this proves the lemma \ref{lm1}.
\hfill$\square$


\noindent 
{\bf Proof of Theorem \ref{thm2}.}
For a given compact interval $I=[a,b]$, consider the function $\tilde{g}$ defined as \eqref{cutoff}. From the above discussion and lemma \ref{lm1},
it follows that $\tilde{g}$ satisfies {\bf (C3)} with the Lipschitz condition 
\begin{eqnarray*}
	{\rm Lip}(\tilde{g})\le\tilde\beta(\omega):= \beta(1+\frac{\max\{|a|, |b|\}}{\omega}).
\end{eqnarray*}
Since $a,b$ are finite, we see that 
$$
\tilde\beta(\omega)\to \beta~~~{\rm as}~~\omega \to +\infty.
$$
Hence,
when $\alpha \ge 2(1-\frac{1}{K})$, by condition \eqref{contracself2}
we can choose a $\epsilon_1 >0$ and a sufficiently large $\omega_1$ 
such that
$\tilde{\beta}(\omega_1) < \beta+\epsilon_1< (K-1)(\alpha K-K+1)$, i.e., $\tilde{\beta}(\omega_1)$ satisfies condition $\eqref{contracself2}$;
when $\alpha< 2(1-\frac{1}{K})$, by condition \eqref{contracselfU1} 
there exists a $\epsilon_2 >0$ and a sufficiently large $\omega_2$ such that 
$\tilde{\beta}(\omega_2) < \beta +\epsilon_2 \le \frac{1}{4} \alpha^2 K^2$, i.e., $\tilde{\beta}(\omega_2)$ satisfies condition \eqref{contracself1}.
Therefore, Theorem \ref{thm1} is available to the functional equation 
\begin{eqnarray}\label{tilde}
\varphi^2(x)=h(\varphi(f(x)))+\tilde{g}(x),
\end{eqnarray}
guaranteeing that 
there exists a bounded continuous function $\tilde{\varphi}$ on $\mathbb{R}$ satisfying \eqref{tilde}.
Restricting equation \eqref{tilde} on $I$, we get that
$\tilde{\varphi}$ satisfies equation $\eqref{rose}$ on $I$.
The proof is completed.
\hfill$\square$


In what follows, we further find continuous solutions of equation $\eqref{rose}$ on the whole $\mathbb{R}$
in the case of unbounded $g$.
We need the following hypotheses:

\begin{description}
	\item  ${(\bf C1^{'})}$ $h$ satisfies condition {(\bf C1)} and there is a real constant $\kappa_h$ such that
	\begin{eqnarray}\label{hlinearbounded}
		\underset{x\in \mathbb{R}}{\sup} |h(x)-\kappa_{h} x| < +\infty;
	\end{eqnarray}

	\item ${(\bf C2^{'})}$ $f$ satisfies condition {(\bf C2)} and there is a real constant $\kappa_f$ such that
	\begin{eqnarray}\label{flinearbounded}
	\underset{x\in \mathbb{R}}{\sup} |f(x)-\kappa_f x| < +\infty;
	\end{eqnarray}

	\item ${(\bf C3^{'})}$ $g$ satisfies Lipschitz condition $\text{Lip}(g)\le \beta$ and 
there is a real constant $\kappa_g\ne 0$ such that
	\begin{eqnarray}\label{glinearbounded}
	\underset{x\in \mathbb{R}}{\sup} |g(x)-\kappa_g x| < +\infty.
	\end{eqnarray}	
\end{description}

For a given constant $\kappa \in \mathbb{R}$, consider
$$
{\cal X}(\mathbb{R};\kappa):=\{\varphi :\mathbb{R} \to \mathbb{R}~|~ \underset{x\in \mathbb{R}}{\sup}|\varphi(x)-\kappa x| < +\infty\},
$$
which is a metric space equipped with
\begin{eqnarray}\label{distance}
	d(\varphi_1,\varphi_2) :=\underset{x\in \mathbb{R}}{\sup} |\varphi_1(x)-\varphi_2(x)|,
	~~~\forall \varphi_1,\varphi_2 \in  {\cal X}(\mathbb{R};\kappa)
\end{eqnarray}
because
$d(\varphi_1,\varphi_2)\le \sup_{x\in \mathbb{R}}|\varphi_1(x)-\kappa x|+\sup_{x\in \mathbb{R}}|\varphi_2(x)-\kappa x| < +\infty$.
Further,
if $\{\varphi_n\}$ is a Cauchy sequence in ${\cal X}(\mathbb{R};\kappa)$, then 
for arbitrary $\epsilon >0$ there exists an integer $N(\epsilon)$
such that 
\begin{eqnarray}\label{Cauchyseq}
	d(\varphi_n,\varphi_m) =\underset{x\in \mathbb{R}}{\sup}|\varphi_n(x)-\varphi_m(x)|< \epsilon,~~~~\forall m,n > N(\epsilon).
\end{eqnarray}
It follows that for each fixed $x\in \mathbb{R}$ the sequence
$\{\varphi_n(x)\}$ is also a Cauchy sequence in $\mathbb{R}$. 
The completeness of $\mathbb{R}$ implies that the limit $\lim_{n \to +\infty}\varphi_n(x)$ exists.
Define $\varphi: \mathbb{R}\to \mathbb{R}$ such that
$$
\varphi(x):=\lim_{n \to +\infty}\varphi_n(x).
$$
First, we claim that 
\begin{eqnarray}
	d(\varphi_n,\varphi) \to 0~~\mbox{ as }~ n\to \infty.
	\label{ddd}
\end{eqnarray}
In fact, it follows from \eqref{Cauchyseq} 
that 
$$
|\varphi_n(x)-\varphi_m(x)|< \epsilon,~~~\forall n,m > N(\epsilon) ,~\forall x\in \mathbb{R}.
$$
Letting $m\to \infty$, one obtains that $|\varphi_n(x)-\varphi(x)|\le \epsilon$ for all $n > N(\epsilon)$
and for all $x\in \mathbb{R}$. Then,
$d(\varphi_n,\varphi)=\underset{x\in \mathbb{R}}{\sup}|\varphi_n(x)-\varphi(x)| \le \epsilon$ for all $n \ge N(\epsilon)$, 
which proves the claim. 
Second, by (\ref{ddd}),
\begin{eqnarray}
	\sup_{x\in \mathbb{R}} |\varphi(x)-\kappa x| 
	&\le& 
	\sup_{x\in \mathbb{R}} \{|\varphi_n(x)-\varphi(x)|+|\varphi_n(x)- \kappa x|\} 
	\nonumber\\
	&\le& 
	\sup_{x\in \mathbb{R}} |\varphi_n(x) -\varphi(x)|+ \sup_{x\in \mathbb{R}} |\varphi_n(x)-\kappa x| < +\infty
	\nonumber
\end{eqnarray}
for sufficiently large $n$, which implies that
$\varphi \in {\cal X}(\mathbb{R};\kappa)$. This proves that
${\cal X}(\mathbb{R};\kappa)$ is complete.

For a constant $L>0$, let
$$
{\cal X}(\mathbb{R};\kappa,L) : = 
{\cal X}(\mathbb{R};\kappa) \cap \{ \varphi : \mathbb{R} \to \mathbb{R}~|~
{\rm Lip}(\varphi) \le L \}.
$$
Clearly, ${\cal X}(\mathbb{R};\kappa,L)$ is a closed subset of 
${\cal X}(\mathbb{R};\kappa)$.


\begin{lm}
	The set ${\cal X}(\mathbb{R};\kappa,L)$ is non-empty if and only if $|\kappa| \le L$.
	\label{lm2}
\end{lm}


\noindent
{\bf Proof}.
Since Lemma \ref{lm2} is clear for $\kappa=0$, we only discuss the case that $\kappa \ne 0$.
Choosing a function  
$
\varphi \in {\cal X}(\mathbb{R};\kappa,L)
$,
we claim that 
Lip$(\varphi) \ge |\kappa|$.
We prove it by contradiction.
Assume that Lip$(\varphi) < |\kappa|$, i.e.,
\begin{eqnarray}
	|\varphi(x) - \varphi(y)|\le {\rm Lip}(\varphi)|x-y| 
	< \frac{|\kappa|+{\rm Lip}(\varphi)}{2} |x-y|,~~~~ \forall x,y \in \mathbb{R}. 
	\label{le}
\end{eqnarray} 
Since 
$
\varphi \in {\cal X}(\mathbb{R};\kappa,L)
$,
we have
$
\varphi(x) = \kappa x + \psi(x)
$,
where
$
\sup_{x\in \mathbb{R}}|\psi(x)| < +\infty
$.
For a monotone $\psi$,
by the boundedness of $\psi$, the limit 
$\lim_{x\to +\infty} \psi(x)$ exists.
Let
\begin{eqnarray}
	A:= \lim_{x\to +\infty} \psi(x) < +\infty.
	\nonumber
\end{eqnarray} 
Then for sufficiently large $x, y(|x-y|=1)$
\begin{eqnarray}
	|\varphi(x) - \varphi(y)| &=& |\kappa x - \kappa y + \psi(x) - \psi(y)| 
	\nonumber\\
	&\ge&
	|\kappa||x-y| - |\psi(x) - \psi(y)|
	\nonumber\\
	&\ge&
	|\kappa||x-y| - \frac{|\kappa|+{\rm Lip}(\varphi)}{2}|x-y|
	\nonumber\\
	&=& \frac{|\kappa|+{\rm Lip}(\varphi)}{2} |x-y|,
	\nonumber
\end{eqnarray}
which contradicts to \eqref{le}.
For a non-monotone $\psi$, we can choose $x,y(x\ne y)$
such that
$\psi(x) = \psi(y)$ and therefore
$$
|\varphi(x) - \varphi(y)| = |\kappa x -\kappa y| = |\kappa||x-y|,
$$
which also contradict to \eqref{le}.
This proves the claim and necessity.
On the other hand,
if 
$
|\kappa| \le L
$, 
obviously, the function 
$
y= \kappa x
$ is contained in 
$
{\cal X}(\mathbb{R};\kappa,L)
$.
This proves the sufficiency and completes the proof of Lemma \ref{lm2}.
\hfill
$\square$


\begin{thm}
Suppose that functions $h,f$ and $g$ satisfy conditions ${(\bf C1^{'})}$, ${(\bf C2^{'})}$ and ${(\bf C3^{'})}$,
where 
constants $K,\alpha$ and $\beta$ satisfy condition \eqref{contracself1} or \eqref{contracself2}.
Then equation $\eqref{rose}$ has a continuous solution on $\mathbb{R}$.
\label{thm3}
\end{thm}


\noindent {\bf Proof}. 
By {\bf (C1)}, $h$ is bijective and the inverse $h^{-1}$ maps $\mathbb{R}$ onto itself.
From the inequality $\eqref{hlinearbounded}$ we get that 
$K\sup_{x\in \mathbb{R}} |\frac{1}{\kappa_h}x-h^{-1}(x)|=
K\sup_{x\in \mathbb{R}} |x-h^{-1}(\kappa_hx)|\le \sup_{x\in \mathbb{R}} |h(x)-h(h^{-1}(\kappa_hx))|=\sup_{x\in \mathbb{R}} |h(x)-\kappa_hx| < +\infty$, i.e., 
\begin{eqnarray}\label{omega_1}
	h^{-1}(x)=\frac{1}{\kappa_h}x+\omega_1(x), 
\end{eqnarray}
where $\omega_1(x)$ is bounded. Similarly, inequalities $\eqref{flinearbounded}$
and $\eqref{glinearbounded}$ imply that
\begin{eqnarray}\label{omega_23}
	f^{-1}(x) =  \frac{1}{\kappa_f}x+\omega_2(x)~~~
	\text{and}~~~
	g(x)  =  \kappa_gx+\omega_3(x),
\end{eqnarray}
where 
$\omega_2(x)$ and $\omega_3(x)$ are bounded.

Define the same mapping $\mathcal{T} :  {\cal X} (\mathbb{R};\kappa,L) \to C^0(\mathbb{R})$ as in $\eqref{TTT}$, i.e.,
$$
\mathcal{T}\varphi(x)= h^{-1}(\varphi^2\circ f^{-1}(x)-g \circ f^{-1}(x)).
$$
By lemma \ref{lm2}, 
in order to ensure that 
$
{\cal X}(\mathbb{R}; \kappa,L)
$
is non-empty, we require that
$$
|\kappa| \le L.
$$

We claim that
\begin{eqnarray}
\sup_{x\in \mathbb{R}}|(\mathcal{T}\varphi)(x) - \kappa x| < +\infty, 
~~~~\forall \varphi \in {\cal X}(\mathbb{R}; \kappa,L)
\nonumber
\end{eqnarray}
if 
\begin{eqnarray}
\kappa = 
\frac{\kappa_h \kappa_f + \sqrt{\kappa_h^2\kappa_f^2 + 4\kappa_g}}{2} 
\left( {\rm or}~
\frac{\kappa_h \kappa_f - \sqrt{\kappa_h^2\kappa_f^2 + 4\kappa_g}}{2} \right)
\label{kappane0}.
\end{eqnarray}
In fact,
every function $\varphi\in {\cal X}(\mathbb{R};\kappa)$ satisfies 
\begin{eqnarray}
\sup_{x\in \mathbb{R}}|\varphi^2(x)-\kappa^2x|
&=&\sup_{x\in \mathbb{R}}|\varphi^2(x)-\kappa\varphi(x) +\kappa\varphi(x)-\kappa^2x| 
\nonumber\\
&\le& \sup_{x\in \mathbb{R}} 
|\varphi^2(x)-\kappa\varphi(x)| +\kappa\sup_{x\in \mathbb{R}}|\varphi(x)-\kappa x|<+\infty.
\nonumber
\end{eqnarray}
It follows from \eqref{omega_1} and \eqref{omega_23}
that 
\begin{eqnarray}
	\begin{split}
		&~\underset{x\in \mathbb{R}}{\sup} |
 \mathcal{T}\varphi (x)
-\kappa x|
		\nonumber\\
		&= \underset{x\in \mathbb{R}}{\sup} |\frac{1}{\kappa_h}(\varphi^2 \circ f^{-1}(x)-g\circ f^{-1}(x))+\omega_1(\varphi^2\circ f^{-1}(x)-g\circ f^{-1}(x))-\kappa  x|
		\nonumber\\
		&=\underset{x\in \mathbb{R}}{\sup} |\frac{1}{\kappa_h}(\kappa^2f^{-1}(x)-g \circ f^{-1}(x))+\frac{1}{\kappa_h}(\varphi^2\circ f^{-1}(x)-\kappa^2f^{-1}(x))
		\nonumber\\
		&~~~+\omega_1(\varphi^2\circ f^{-1}(x)-g\circ f^{-1}(x))-\kappa  x|
		\nonumber\\
		&=\underset{x\in \mathbb{R}}{\sup} |\frac{1}{\kappa_h}(\kappa^2f^{-1}(x)-\kappa_gf^{-1}(x))+\frac{1}{\kappa_h}(\varphi^2\circ f^{-1}(x)-\kappa^2f^{-1}(x))-\frac{1}{\kappa_h}\omega_3(f^{-1}(x))
		\nonumber\\
		&~~~ +\omega_1(\varphi^2\circ f^{-1}(x)-g\circ f^{-1}(x))-\kappa x|
		\nonumber\\
		&=\underset{x\in \mathbb{R}}{\sup} |(\frac{\kappa^2-\kappa_g}{\kappa_h\kappa_f} -\kappa)x +\frac{\kappa^2-\kappa_g}{\kappa_h}\omega_2(x) +\frac{1}{\kappa_h}(\varphi^2\circ  f^{-1}(x)-\kappa^2f^{-1}(x))
		\nonumber\\
		&~~~ -\frac{1}{\kappa_h}\omega_3(f^{-1}(x))+\omega_1(\varphi^2\circ f^{-1}(x)-g\circ f^{-1}(x))|< +\infty
	\end{split}
	\end{eqnarray}
if $\kappa$ is chosen such that 
\begin{eqnarray}\label{bounded condition}
\frac{\kappa^2-\kappa_g}{\kappa_h\kappa_f}-\kappa =0.
\end{eqnarray}
By lemma \ref{lm2}, we obtain that
\begin{eqnarray}
K \le |\kappa_h|,~~ \alpha \le |\kappa_f|,
~~ |\kappa_g| \le \beta,
\label{Lipkappa}
\end{eqnarray}
which implies 
\begin{eqnarray}
\kappa_h^2\kappa_f^2 + 4\kappa_g
\ge
\alpha^2K^2 - 4\beta. 
\label{ge4}
\end{eqnarray}
Under condition \eqref{contracself1} or 
\eqref{contracself2}, as shown in the proof of Theorem \ref{thm1}, condition \eqref{Delta} holds, i.e.,
$
\alpha^2K^2 - 4\beta \ge 0
$
and therefore by \eqref{ge4} 
$
\kappa_h^2\kappa_f^2 + 4\kappa_g \ge 0
$.
Hence, two real roots of equation \eqref{bounded condition} are
\begin{eqnarray}
\kappa_{1,2} = \frac{\kappa_h \kappa_f \pm \sqrt{\kappa_h^2\kappa_f^2 + 4\kappa_g}}{2},
\nonumber
\end{eqnarray}
i.e., condition \eqref{kappane0}.
This completes the proof of the claim.

Further,
as proved in Theorem \ref{thm1}, 
\begin{eqnarray}
{\rm Lip}(\mathcal{T}\varphi) \le L,
~~~~ \forall \varphi \in {\cal X}(\mathbb{R}; \kappa,L)
\nonumber
\end{eqnarray}
if $L$ satisfies \eqref{self}, and 
\begin{eqnarray}
{\rm Lip}(\mathcal{T}) < 1
\nonumber
\end{eqnarray}
if $L$ satisfies \eqref{contrac}.

Consequently, $\mathcal{T}$ is a contractive self-mapping on
a non-empty set
$
{\cal X}(\mathbb{R}; \kappa_*,L),
$
where 
$\kappa_* = \kappa_1$ as $|\kappa_1| \le |\kappa_2|$;
$\kappa_* = \kappa_2$ as $|\kappa_1| > |\kappa_2|$,
if 
$L$ is chosen to fulfill \eqref{self}, \eqref{contrac}
and 
\begin{eqnarray}
|\kappa_*| \le L.
\label{nonemp2}
\end{eqnarray}
We claim that such $L$ exists. In fact,
under condition \eqref{contracself1} or \eqref{contracself2}, 
from the proof of Theorem \ref{thm1} we can see that
\begin{eqnarray}
\frac{1}{2}\alpha K-\frac{1}{2}\sqrt{\alpha^2 K^2-4\beta} < K-1,
\nonumber
\end{eqnarray}
as shown in \eqref{intersection}.
By \eqref{Lipkappa} and \eqref{ge4}, we see that
\begin{eqnarray}
|\kappa_*|& =& \left|\frac{\kappa_h \kappa_f - \sqrt{\kappa_h^2\kappa_f^2 + 4\kappa_g}}{2} \right|
=\left|\frac{2\kappa_g}{\kappa_h \kappa_f + \sqrt{\kappa_h^2\kappa_f^2 + 4\kappa_g}} \right|
\nonumber\\
&\le&
\frac{2\beta}{\alpha K + \sqrt{\alpha^2K^2 - 4 \beta}} 
\nonumber\\
&=&
\frac{\alpha K - \sqrt{\alpha^2K^2 - 4 \beta }}
{2} ~~~~~~ {\rm as}~~\kappa_h \kappa_f >0,
\label{le2}
\\
|\kappa_*|& =& \left|\frac{\kappa_h \kappa_f + \sqrt{\kappa_h^2\kappa_f^2 + 4\kappa_g}}{2} \right|
=\left|\frac{2\kappa_g}{\kappa_h \kappa_f - \sqrt{\kappa_h^2\kappa_f^2 + 4\kappa_g}} \right|
\nonumber\\
&\le&
\frac{2\beta}{|\kappa_h \kappa_f| + \sqrt{\kappa_h^2\kappa_f^2 + 4\kappa_g}}
\le 
\frac{2\beta}{\alpha K + \sqrt{\alpha^2K^2 - 4 \beta }} 
\nonumber\\
&=& 
\frac{\alpha K - \sqrt{\alpha^2K^2 - 4 \beta }}
{2}~~~~~~~{\rm as}~~ \kappa_h \kappa_f \le 0.
\label{le3}
\end{eqnarray}
Thus, from \eqref{intersection}, \eqref{le2}
and 
\eqref{le3} we can choose a $L$ satisfying
\eqref{self}, \eqref{contrac} and \eqref{nonemp2}.
This proves the claim.
It follows that there exist constants $\kappa$ and $L$ such that $\mathcal{T}$ is a contractive self-mapping on the non-empty set ${\cal X}(\mathbb{R};\kappa,L)$,
implying that a unique solution of $\eqref{rose}$ exists in ${\cal X}(\mathbb{R};\kappa,L)$. This completes the proof. 
\hfill$\square$


The continuous solution found in Theorem~\ref{thm3} is unbounded.
Otherwise, $\varphi^2$ is bounded and satisfies
$
	\varphi^2(x)=h(\varphi(f(x)))+g(x),
$
but the right hand side is unbounded because $h$ is continuous yields that $h(\varphi(f(x)))$ is bounded and 
condition ${(\bf C3^{'})}$ implies that $g$ is unbounded. 
This is a contradiction.

\begin{exam}{\rm
	Theorem \ref{thm3} can be applied to the following equation
	$$
	\varphi^2(x)= -2\varphi(2x)+x+\sin x,~~~~ x\in \mathbb{R},
	$$
which is of the form (\ref{rose}),
where $h(x)=-2x, f(x)=2x$ and $g(x)=x+\sin x$.
One can verify
conditions ${(\bf C1^{'})}$-${(\bf C3^{'})}$ with  
$K=\alpha=\beta=\kappa_h=\kappa_f=2$ and $\kappa_g=1$. Further, 
$$
	-\frac{1}{4}(\kappa_f\kappa_h)^2=-4 < 1.
$$
It is the same as in Example $1$ that constants $K, \alpha$ and $\beta$ satisfy \eqref{contracself2}.
By Theorem \ref{thm3}, the equation has a continuous solution on $\mathbb{R}$.
}\label{exam2}
\end{exam}


\section{Case without Lipschitz conditions
}

In this section we consider the case where we do not impose the Lipschitz condition to $g$ and the inverses of $h$ and $f$.
In this case we hardly use a fixed point theorem, but more 
solutions of equation \eqref{rose} can be constructed piecewise as follows.


The following theorem is devoted to the increasing case, that is,
functions $h,f$ and $g$ are all strictly increasing and continuous.
\begin{figure}
\centering
\includegraphics[width=2.2in]{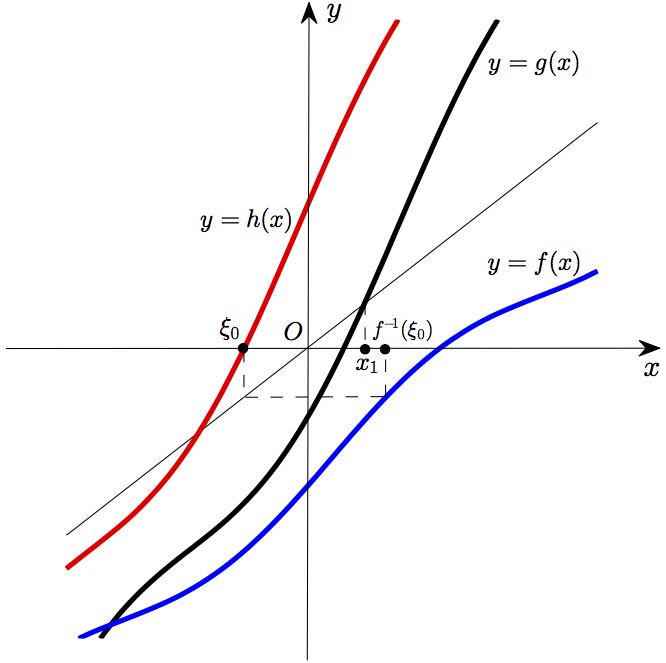}
\vspace{1pt}
\caption{Graphs of $f,g$ and $h$.}
\label{fgh}
\end{figure}

\begin{thm}
Suppose that functions $h,f$ and $g$ are all strictly increasing and continuous on $\mathbb{R}$,
$h: \mathbb{R}\to \mathbb{R}$ is surjective, $f(x) < x$ for all $x\in \mathbb{R}$, and 
$g$ has a fixed point $x_1$ such that 
\begin{eqnarray}\label{x1}
\xi_0< x_1 \le f^{-1}(\xi_0),
\end{eqnarray}
where $\xi_0$ is the unique zero of $h$,
and $g(x) \ge x$ as $x \ge x_1$.
Then any two  strictly increasing and continuous surjections $\varphi_0 : [x_0, x_1] \to [x_1,x_2]$ 
and $\varphi_1 : [x_1,x_2] \to [x_2,x_3]$,
where 
$x_0:=f(x_1)$, 
$x_3:=h(x_1)+ x_1$,
and
$x_2$ is chosen arbitrarily such that 
\begin{eqnarray}\label{x2}
x_1 < x_2 < h(x_1) + x_1,
\end{eqnarray}
can be extended uniquely to a continuous solution of $\eqref{rose}$ on $\mathbb{R}$.
\label{thm4}
\end{thm}

The relationship among $f, g$ and $h$ required in Theorem \ref{thm4} can be shown intuitively in Figure \ref{fgh}.
It is easy to find such functions $f, g$ and $h$, for example, $f(x)=x-1$, $g(x)=2x$ and $h(x)=x+{1}/{2}$. 
Clearly, they are all strictly increasing and continuous, 
$h(\mathbb{R})=\mathbb{R}$, $f(x) < x$ for all $x\in \mathbb{R}$, and $g$ has a fixed point $x_1=0$, i.e., $g(0)=0$.
Moreover, $h$ has a unique zero $\xi_0=-1/2$. One can check that
$f^{-1}(\xi_0)=1/2>x_1>\xi_0$ and that
$g(x)=2x \ge x$ for all $x\ge 0$.
Hence, $f, g$ and $h$ 
satisfy the conditions of Theorem \ref{thm4}.


\noindent
{\bf Proof of Theorem~\ref{thm4}.} 
First, we construct a solution of \eqref{rose} on $[x_0, +\infty)$.
Since $f(x) < x$ for all $x\in \mathbb{R}$, it is clear that $x_0 = f(x_1) < x_1$. 
Because $h$ is strictly increasing, we see from condition \eqref{x1} that $x_1 = h(\xi_0)+x_1 < h(x_1)+x_1=x_3$,
which implies that the choice of $x_2$ in \eqref{x2} is reasonable.

Having given $x_0< x_1 <x_2 <x_3$ above, we know that there are infinitely many 
increasing homeomorphisms $\varphi_0: [x_0, x_1] \to [x_1,x_2]$
and $\varphi_1 : [x_1,x_2] \to [x_2,x_3]$.
Define
\begin{eqnarray}\label{varphi2}
	\varphi_2(x) := h\circ \tilde{\varphi}_1\circ f\circ \varphi^{-1}_1(x) + g\circ \varphi_1^{-1}(x), ~~x\in [x_2,x_3],
\end{eqnarray}
where 
\begin{eqnarray}
\tilde{\varphi}_1(x):=\begin{cases}
\varphi_0(x),~~~~        & x\in [x_0, x_1),
\\
\varphi_1(x),~~~~        & x\in [x_1,x_2].
\end{cases}
\label{def1}
\end{eqnarray}
Clearly, $\tilde{\varphi}_1$ is well defined and strictly increasing continuous.
By the choice of $\varphi_1$ and the assumption of $f$, we get that $f \circ \varphi_1^{-1}(x) \in [x_0,x_2]$ as $x\in [x_2,x_3]$, 
i.e., $f \circ \varphi_1^{-1}([x_2,x_3])$ is a subset of the domain of $\tilde{\varphi}_1$, implying that $\varphi_2$ is well defined.
Let
\begin{eqnarray}\label{x4x4}
x_4:=\varphi_2(x_3).
\end{eqnarray} 
 Then $\varphi_2 : [x_2,x_3] \to [x_3,x_4]$ is an increasing homeomorphism.
 In fact, $\varphi_2$ is strictly increasing continuous because so are all functions on the right hand side of \eqref{varphi2}.
Moreover,  $\varphi_2$ is surjective because we have \eqref{x4x4} and
\begin{eqnarray*}\label{x3}
\varphi_2(x_2) =h\circ \tilde{\varphi}_1\circ f(x_1) +g(x_1)=h\circ \varphi_0(x_0)+x_1=h(x_1)+x_1=x_3,
\end{eqnarray*}
which is obtained from \eqref{def1} and the fact that $x_1$ is a fixed point of $g$.

Assume that for integer $k\ge 2$ a strictly increasing sequence $\{x_i\}_{i=0}^{k+2}$
and $k-1$ increasing homeomorphisms 
$\varphi_i: [x_i,x_{i+1}] \to [x_{i+1},x_{i+2}], i=2,3,\cdots,k$, are well defined 
 such that
\begin{eqnarray}\label{varphi-i}
\varphi_i(x)= h\circ \tilde{\varphi}_{i-1} \circ f \circ \varphi^{-1}_{i-1}(x) + g\circ \varphi_{i-1}^{-1}(x),~~~x\in [x_i,x_{i+1}],
\end{eqnarray}
where
\begin{eqnarray*}\label{tilde-i-1}
\tilde{\varphi}_{i-1}(x):=\begin{cases}
\varphi_0(x), &~~~~ x\in [x_0,x_1),
\\
\varphi_1(x), &~~~~ x\in [x_1,x_2),
\\
\cdots
\\
\varphi_{i-1}(x), &~~~~ x\in [x_{i-1},x_i].
\end{cases}
\end{eqnarray*}
Define
\begin{eqnarray}\label{induction}
\varphi_{k+1}(x):=h\circ \tilde{\varphi}_k \circ f \circ \varphi_k^{-1}(x)+g\circ \varphi_k^{-1}(x), ~~ x\in [x_{k+1},x_{k+2}],
\end{eqnarray}
where
\begin{eqnarray*}
\tilde{\varphi}_k(x):=\begin{cases}
\varphi_0(x), &~~~~ x\in [x_0,x_1),
\\
\varphi_1(x), &~~~~ x\in [x_1,x_2),
\\
\cdots
\\
\varphi_k(x), & ~~~~  x\in [x_k,x_{k+1}].
\end{cases}
\end{eqnarray*}
Obviously, 
$\tilde{\varphi}_k$ is well defined and strictly increasing continuous because $\varphi_i$ is an increasing homeomorphism
 from $[x_i,x_{i+1}]$ to $[x_{i+1},x_{i+2}]$ for each $0\le i \le k$.
By the inductive assumption, we see that the inverse $\varphi_k^{-1}$ is well defined on $[x_{k+1},x_{k+2}]$ 
and $\varphi_k^{-1}([x_{k+1},x_{k+2}]) = [x_k,x_{k+1}]$.
Since $f$ is strictly increasing and satisfies $f(x) < x$ for all $x\in \mathbb{R}$,
we have 
$$
x_0 =f(x_1) < f(x_k) \le f\circ \varphi_k^{-1}(x) \le f(x_{k+1}) < x_{k+1}~~ \text{as}~~ x\in [x_{k+1},x_{k+2}],
$$
i.e., $f \circ \varphi_k^{-1}([x_{k+1},x_{k+2}])$ is a subset of the domain of $\tilde{\varphi}_k$, implying that
$\varphi_{k+1}$ is well defined. 
Letting 
\begin{eqnarray}\label{x_k+3}
x_{k+3}:= \varphi_{k+1}(x_{k+2}),
\end{eqnarray}
we claim that $\varphi_{k+1} : [x_{k+1}, x_{k+2}] \to [x_{k+2}, x_{k+3}]$ is an
increasing homeomorphism. In fact, 
$\varphi_{k+1}$ is strictly increasing continuous because all functions on the right hand side of \eqref{induction} are strictly increasing continuous.
Moreover, $\varphi_{k+1}$ is surjective because we have \eqref{x_k+3} and 
\begin{eqnarray} 
\varphi_{k+1}(x_{k+1})  
&=& 
h\circ \tilde{\varphi}_k \circ f \circ \varphi_k^{-1}(x_{k+1})+g\circ \varphi_k^{-1}(x_{k+1})
\nonumber\\
&=&  h\circ \tilde{\varphi}_k\circ f(x_k)+g(x_k) 
\nonumber\\
&=&  h\circ \tilde{\varphi}_{k-1}\circ f(x_k)+g(x_k)
\nonumber\\
&=& h\circ \tilde{\varphi}_{k-1}\circ f \circ \varphi_{k-1}^{-1}(x_{k+1})+g \circ \varphi_{k-1}^{-1}(x_{k+1})
\nonumber\\
&=&\varphi_k(x_{k+1})=x_{k+2},
\label{varphik1}
\end{eqnarray}
which is deduced from \eqref{induction} and the inductive assumption.
This proves the claim.
Hence, we have proved by induction that
there is 
a strictly increasing sequence 
$\{x_i\}_{i= 0}^{+\infty}$ and a sequence of functions $\{\varphi_i\}_{i\ge 0}$, 
where 
$\varphi_i : [x_i,x_{i+1}] \to [x_{i+1},x_{i+2}]$,
defined by \eqref{varphi-i}, 
is an increasing homeomorphism for each $i\ge 0$.

We further claim that 
\begin{eqnarray}
x_k \to +\infty ~~\mbox{ as } ~~k \to +\infty. 
\label{limit-infty}
\end{eqnarray}
If it is not true, 
let $x_k \to x_*$ as $k \to +\infty$ by the monotonicity.
Putting $x=x_{k+1}$ in \eqref{induction}, we get
\begin{eqnarray}
\label{induction*}
\varphi_{k+1}(x_{k+1})=h\circ \tilde{\varphi}_k \circ f \circ \varphi_k^{-1}(x_{k+1})+g\circ \varphi_k^{-1}(x_{k+1}), 
\end{eqnarray}
where 
$\varphi_{k+1}(x_{k+1})=x_{k+2}$ and $ \varphi_k^{-1}(x_{k+1})=x_k$. 
Since $x_0=f(x_1)<f \circ \varphi_k^{-1}(x_{k+1})=f(x_k)<x_k$, we have
$x_1=\tilde{\varphi}_k(x_0)<\tilde{\varphi}_k \circ f \circ \varphi_k^{-1}(x_{k+1})<\tilde{\varphi}_k(x_k)=x_{k+1}$.
It follows by the strictly increasing monotonicity that 
$\tilde{\varphi}_k \circ f \circ \varphi_k^{-1}(x_{k+1}) \to \tilde{x}$ 
as $k\to+\infty$, where $\tilde{x}\in (x_1,  x_*]$.
Letting $k\to+\infty$ in (\ref{induction*}),
by continuity we obtain
\begin{eqnarray}\label{x_*}
x_* = h(\tilde{x}) + g(x_*).
\end{eqnarray}
On the other hand, 
$h$ is strictly increasing, $h(x_1) > 0$, and 
$g(x) \ge x$ as $x\ge x_1$, 
which imply that $h(\tilde{x}) + g(x_*) > h(x_1)+x_* > x_*$, a contradiction to \eqref{x_*}.
The claimed (\ref{limit-infty}) implies that 
$$
[x_0, +\infty)=\cup_{i=0}^\infty [x_i,x_{i+1}).
$$
Then, define 
\begin{eqnarray}
\varphi_{*}(x) := \varphi_i(x), ~~\forall x\in [x_i,x_{i+1}),~~~ i\ge 0.
\label{phi-def}
\end{eqnarray}
The above discussion shows that the function $\varphi_*$ is well defined and strictly increasing continuous on $[x_0,+\infty)$.
Furthermore, for an arbitrary $x\in [x_1,+\infty)$, there exists an integer $i\ge 2$ such that
$x\in [x_{i-1},x_i)$.
By the definition \eqref{varphi-i} of $\varphi_i$ and the definition \eqref{tilde-i-1} of $\tilde{\varphi}_{i-1}$, 
\begin{eqnarray}
		{\varphi}_*^2(x)\!\!\!\!\!&& =\varphi_i\circ \varphi_{i-1}(x) 
		\nonumber\\
		&&= h\circ \tilde{\varphi}_{i-1} \circ f \circ \varphi^{-1}_{i-1} (\varphi_{i-1}(x)) + g\circ \varphi_{i-1}^{-1}(\varphi_{i-1}(x))
		\nonumber\\
		&&= h\circ \tilde{\varphi}_{i-1} \circ f (x) + g(x)
		\nonumber\\
		&&=h\circ \varphi_* \circ f(x)+g(x),
		\label{check*}
\end{eqnarray} 
implying that function $\varphi_*$ is a solution of equation \eqref{rose} on $[x_1, +\infty)$.


Next, we extend the solution $\varphi_{*}$ from $[x_0,+\infty)$ to the whole real line $(-\infty, +\infty)$. 
Let $x_{-i}:=f^i(x_0), i=1,2,3,\cdots.$ Then
the sequence $\{x_{-i}\}_{i\ge 1}$ is strictly decreasing
and satisfies $x_{-i} \to -\infty$ as $i \to +\infty$ since $f(x) < x$ for all $x\in \mathbb{R}$. 
It gives the partition
$$
(-\infty, x_0)=\cup_{i=0}^\infty [x_{-i-1}, x_{-i}).
$$
For each integer $k\ge 1$ define 
\begin{eqnarray}\label{iii}
\varphi_{-k}(x):= h^{-1}(\varphi_* \circ \varphi_{-k+1} \circ f^{-1} (x)  - g\circ f^{-1}(x) ), ~~x\in [x_{-k}, x_{-k+1}],
\end{eqnarray}
recursively with $\varphi_0$ being $\varphi_*$ on $[x_0, x_1]$,
where $\varphi_*$ is the solution on $[x_0,+\infty)$ obtained in \eqref{phi-def}.
We claim that every $\varphi_{-k}$ is well defined and continuous on $[x_{-k}, x_{-k+1}]$ such that
\begin{eqnarray}
&& \varphi_{-k}(x) > \xi_0,~ \forall x\in [x_{-k}, x_{-k+1}],
\label{returnback}
\\
&& \varphi_{-k}(x_{-k+1}) = \varphi_{-k+1}(x_{-k+1}).
\label{conjunction}
\end{eqnarray}

In fact, for $k=1$ we can see that $\varphi_{-1}$, defined by
\begin{eqnarray}\label{-1}
\varphi_{-1}(x):= h^{-1}(\varphi_*^2 \circ f^{-1} (x)  - g\circ f^{-1}(x) ), ~~x\in [x_{-1}, x_0],
\end{eqnarray}
as in (\ref{iii})
 is well defined  because 
$
f^{-1}([x_{-1}, x_0]) = [x_0, x_1] \subset [x_0, +\infty),
$
i.e., $f^{-1}([x_{-1}, x_0])$ is contained in the domain of $\varphi_*$.
The continuity of $\varphi_{-1}$ comes from the fact that functions on the right hand side of \eqref{-1} are all continuous. 
In order to prove (\ref{returnback}) with the index $-1$ in place of $-k$, we note that 
\begin{eqnarray*}
&&\varphi_*^2 \circ f^{-1}(x)\ge \varphi_*^2 \circ f^{-1}(x_{-1})=x_2,
~~~~~~~~~ \forall x\in [x_{-1}, x_0],
\\
&&g\circ f^{-1}(x) \le g\circ f^{-1}(x_0)= g(x_1) =x_1,
~~~\forall x\in [x_{-1}, x_0],
\end{eqnarray*}
since functions $\varphi_*, f^{-1}$ and $g$ are all strictly increasing and $g(x_1) =x_1$.
It follows from (\ref{-1}) that
\begin{eqnarray}
\varphi_{-1}(x)\!\!\!\!\! && = h^{-1}(\varphi_*^2 \circ f^{-1} (x)  - g\circ f^{-1}(x) )
\nonumber\\
&& \ge h^{-1}(x_2 -x_1) 
\nonumber\\
&&> h^{-1}(0) = \xi_0,
~~~~~~
\forall x\in [x_{-1}, x_0],
\label{gxi0}
\end{eqnarray}
by the definition of $\xi_0$ and the monotonicity of $h$.
This proves (\ref{returnback}) for $k=1$.
Further, from (\ref{-1}) we have
\begin{eqnarray*}
\begin{split}
\varphi_{-1}(x_0) & = h^{-1}(\varphi_*^2 \circ f^{-1} (x_0)  - g\circ f^{-1}(x_0) )
\\
& = h^{-1}( \varphi_1 \circ \varphi_0 (x_1)  - g(x_1) )
\\
& = h^{-1}( x_3  - x_1 )
\\
& =h^{-1} \circ h( x_1 )= x_1 =\varphi_*(x_0)
\end{split}
\end{eqnarray*}
by the choice of $x_3$, 
which proves (\ref{conjunction}) for $k=1$.

Generally assume that
for an integer $k\ge 1$ function $\varphi_{-k}$
is well defined by \eqref{iii}
and continuous on $[x_{-k}, x_{-k+1}]$ such that \eqref{returnback} and \eqref{conjunction}.
Let
\begin{eqnarray}\label{-k-1}
\varphi_{-k-1}(x):= h^{-1}(\varphi_* \circ \varphi_{-k} \circ f^{-1} (x)  - g\circ f^{-1}(x) ), \,~x\in [x_{-k-1}, x_{-k}],
\end{eqnarray}
where $\varphi_*$ is obtained in \eqref{phi-def}.
By \eqref{returnback} and \eqref{x1} we see that
\begin{eqnarray}\label{xi0}
\varphi_{-k} \circ f^{-1} (x)  > \xi_0 \ge f(x_1) =x_0,~~~~ \forall x\in  [x_{-k-1}, x_{-k}],
\end{eqnarray}
i.e., $\varphi_{-k} \circ f^{-1} ([x_{-k-1}, x_{-k}]) $ is contained in the domain of $\varphi_*$,
which implies that $\varphi_{-k-1}$ is well defined.
$\varphi_{-k-1}$ is continuous because all functions on the
right hand side of \eqref{-k-1} are continuous.
Note that $g\circ f^{-1}(x) \le x_1$ for all $x\in [x_{-k-1}, x_{-k}]$ because $g(x_1)=x_1$ and $g$ is strictly increasing. 
It follows from \eqref{-k-1} and \eqref{xi0} that 
\begin{eqnarray*}
	\begin{split}
		\varphi_{-k-1}(x) & = h^{-1}(\varphi_* \circ \varphi_{-k} \circ f^{-1} (x)  - g\circ f^{-1}(x) )
		\\
		& > h^{-1}(\varphi_*(x_0) - x_1) = h^{-1}(0)= \xi_0,~~~~~ \forall x\in [x_{-k-1}, x_{-k}],
	\end{split}
\end{eqnarray*}
by the monotonicity of functions $h$ and $\varphi_*$. This proves \eqref{returnback} for the index $-k-1$.
Furthermore,
by \eqref{-k-1}, \eqref{conjunction} and the definition \eqref{iii} of $\varphi_{-k}$, 
we obtain
\begin{eqnarray*}\label{-k+1}
\begin{split}
\varphi_{-k-1}(x_{-k}) & = h^{-1}(\varphi_* \circ \varphi_{-k} \circ f^{-1} (x_{-k})  - g\circ f^{-1}(x_{-k}) )
\\
& = h^{-1}( \varphi_* \circ \varphi_{-k} (x_{-k+1})  - g(x_{-k+1}) )
\\
& = h^{-1}( \varphi_* \circ \varphi_{-k+1}(x_{-k+1})  - g(x_{-k+1}) )
\\
& = h^{-1}( \varphi_* \circ \varphi_{-k+1} \circ f^{-1}(x_{-k})  - g \circ f^{-1}(x_{-k}) )
\\
& =\varphi_{-k}(x_{-k}),
\end{split}
\end{eqnarray*}
which proves \eqref{conjunction} for the index $-k-1$ and completes the proof of the claim.


Finally, define a function $\varphi$ on $\mathbb{R}$ by 
$$
\varphi(x) := \varphi_i(x),~~~~\forall x\in [x_i,x_{i+1}),~~ i \in \mathbb{Z}. 
$$
Then, $\varphi$ is continuous on $\mathbb{R}$ by $\eqref{conjunction}$ because 
$\varphi(x)=\varphi_*(x)$ for all $x\in [x_0, +\infty)$, as defined in \eqref{phi-def}.
We have checked that $\varphi$ satisfies equation \eqref{rose} for all $x \in [x_1, +\infty)$ in \eqref{check*}.
For an arbitrary $x \in (-\infty, x_1)$, without loss of generality, $x\in [x_{-k+1},x_{-k+2})$ for a certain integer $k\ge 1$,
by \eqref{returnback} and \eqref{iii} we have
\begin{eqnarray*}
{\varphi}^2(x)  =\varphi_* \circ \varphi_{-k+1}(x) = h \circ \varphi_{-k} \circ f (x)  + g(x)=h\circ \varphi \circ f(x)+g(x),
\end{eqnarray*}
i.e., function $\varphi$ satisfies equation \eqref{rose} for all $x \in (-\infty, x_1)$. 
It follows that $\varphi$ is a continuous solution of \eqref{rose} on $\mathbb{R}$.


In order to prove the uniqueness of $\varphi$, assume that another function $\hat{\varphi}$, 
which is defined on $\mathbb{R}$ and coincides with $\varphi_0$ and $\varphi_1$
on $[x_0,x_1]$ and  $[x_1,x_2]$ respectively, 
also satisfies  
equation \eqref{rose} for all $x\in \mathbb{R}$. 
Restricting equation \eqref{rose} to the interval $[x_1,x_2]$, 
we obtain that 
$$
\hat{\varphi} \circ \varphi_1(x)=h \circ \tilde{\varphi}_1 \circ f(x) + g(x)~~\text{as}~~ x\in [x_1,x_2],
$$
or equivalently say,
$$
\hat{\varphi} (x)=h \circ \tilde{\varphi}_1 \circ f \circ \varphi_1^{-1}(x) + g \circ \varphi_1^{-1}(x) ~~\text{as}~~
 x\in [x_2,x_3].
$$
It follows from \eqref{varphi2} that $\hat{\varphi}|_{[x_2,x_3]}\equiv \varphi_2$. Further,
by induction we can prove that 
\begin{eqnarray}
\hat{\varphi}|_{[x_i,x_{i+1}]}\equiv \varphi_i,~~~\forall i\ge 3.
\label{+++}
\end{eqnarray}
On the other hand, restricting equation \eqref{rose} to the interval $[x_0,x_1]$, 
we obtain
$$
\varphi_1 \circ \varphi_0(x)=h \circ \hat{\varphi} \circ f(x) + g(x)~~\text{as}~~ x\in [x_0,x_1],
$$
or equivalently say,
$$
\hat{\varphi}(x)= h^{-1}(\varphi_*^2 \circ f^{-1} (x)  - g\circ f^{-1}(x) )
~~\text{as}~~ x\in [x_{-1}, x_0].
$$
By \eqref{-1} we get 
$\hat{\varphi}|_{[x_{-1}, x_0]}\equiv \varphi_{-1}$. By induction one can prove that 
\begin{eqnarray}
\hat{\varphi}|_{[x_{-i},x_{-i+1}]}\equiv \varphi_{-i},~~~\forall i\ge 2. 
\label{---}
\end{eqnarray}
It follows from (\ref{+++}) and (\ref{---}) that $\hat{\varphi}\equiv \varphi$,
implying the uniqueness of $\varphi$. This completes the proof.
\hfill$\square$

Theorem \ref{thm4} has some overlaps with Theorem~\ref{thm3}.
Like Theorem~\ref{thm3}, it also deals with unbounded $g$ 
since it requires $g(x)\ge x$ for all $x\ge x_1$.
Theorems \ref{thm3} and \ref{thm4} are both applicable to given functions $h(x)=3x+1, f(x)=x-1$ and $g(x)=x$, but 
Theorem \ref{thm4} gives more solutions.
However,
Theorem \ref{thm4} can be applied to functions $h(x)=x+{1}/{2}$, $f(x)=x-1$ and $g(x)=2x$, as illustrated 
just below Theorem~\ref{thm4},
but Theorem \ref{thm3} cannot
because Theorem \ref{thm3} requires $h$ to be expansive.
This does not mean that the conditions of Theorem \ref{thm4} are weaker.
For example,
Theorem \ref{thm4}
can not be applied to the given functions $h(x)=-2x, f(x)=2x$ and $g(x)=x+\sin x$, 
which were considered with Theorem \ref{thm3} in Example \ref{exam2},
because Theorem \ref{thm4} requires that $h$ is strictly increasing.


\section{ Some Remarks}

In the proof of Theorem~\ref{thm4} we used two methods in construction of solutions.
One is the usual method of ``first locate points then define functions'' as used on $(-\infty, x_0)$.
The other is the method of ``locate point and define function alternately'' as done on $[x_0,+\infty)$.
If we use the method of ``first locate points then define functions'' on $[x_0,+\infty)$ and,
similarly to our construction on $(-\infty, x_0)$,
locate
$$
x_i:=f^{-i}(x_0),~~~\forall i \ge 1,
$$ 
we have the partition
$[x_0,+\infty)=\cup_{i=0}^{+\infty}[x_i,x_{i+1})$,
provided that $f$ is a homeomorphism additionally.
In the routine of construction, 
for arbitrarily chosen strictly increasing homeomorphisms
$\varphi_0: [x_0,x_1] \to [x_1, x_2]$ 
and 
$\varphi_1 : [x_1,x_2] \to [x_2,x_3]$, we define 
 \begin{eqnarray}\label{k+1}
 \varphi_{i}(x):= h \circ \varphi_{i-2} \circ f \circ \varphi_{i-1}^{-1} (x) + g \circ \varphi_{i-1}^{-1} (x),~~~ \forall x\in [x_{i}, x_{i+1}],
 \end{eqnarray}
for all integers $i\ge 2$ inductively and connect them to make a continuous solution.
We can prove that $\varphi_{i}: [x_{i}, x_{i+1}] \to [x_{i+1}, x_{i+2}]$
is an increasing homeomorphism if 
\begin{eqnarray*}
h(x_{i-1}) + g(x_{i-1}) = x_{i+1}~~~ {\rm and} ~~~ h(x_{i}) + g(x_{i}) = x_{i+2},
\end{eqnarray*}
which actually impose a strong condition on $h$ and $g$ at each point of the sequence $\{x_i\}_{i\ge 1}$.


Theorem \ref{thm4} requires two conditions: the fixed point $x_1$ of $g$ is chosen to fulfill \eqref{x1}, i.e., $\xi_0 < x_1 \le f^{-1}(\xi_0)$,
and
\begin{eqnarray}
g(x) \ge x~~~\mbox{as}~~x\ge x_1.
\label{gxg}
\end{eqnarray}
{\it If we do not consider \eqref{x1}, the existence of continuous solutions of equation \eqref{rose} 
remains unknown.}
Actually,
if $x_1 \le \xi_0$, then $h(x_1) \le 0$, i.e., 
the inequality $x_1 < h(x_1) +x_1$ is not true,
which implies that
there does not exist $x_2$ satisfying \eqref{x2}
and therefore our construction, which depends on \eqref{x2} because we require $x_1<x_2<x_3=\varphi_2(x_2)=h(x_1)+x_1$ as shown in \eqref{x3},
does not work.
On the other hand,
if $x_1 > f^{-1}(\xi_0)$, we have $x_0=f(x_1) > f(f^{-1}(\xi_0)) =\xi_0$.
Then we cannot use the same method of construction as in Theorem \ref{thm4} on $(-\infty, x_0)$
because, when defining $\varphi_{-2}$, as doing in \eqref{iii} with $k=2$, we need
\begin{eqnarray*}
 \varphi_{-1} \circ f^{-1}(x) \ge x_0, ~~~~\forall x\in [x_{-2},x_{-1}],
\end{eqnarray*}
which however is not guaranteed by the inequality 
$\xi_0 < x_0$ and \eqref{gxi0} (i.e., $\varphi_{-1}(x) > \xi_0$ for all $x\in [x_{-1},x_0]$).
Moreover, as $x_{-1} \le \xi_0 < x_0$, we can define $\varphi_{-2}$ by
\begin{eqnarray}\label{-2}
\varphi_{-2}(x):=h^{-1}(\tilde{\varphi}_{-1}^2 \circ f^{-1}(x)-g\circ f^{-1}(x)),~~~~ \forall x\in [x_{-2}, x_{-1}],
\end{eqnarray}
where 
\begin{eqnarray*}
\tilde{\varphi}_{-1}(x):= \begin{cases}
\varphi_{-1}(x),~~~~  & x\in [x_{-1}, x_0),
\\
\varphi_*(x),~~~~  & x\in [x_0, +\infty).
\end{cases}
\end{eqnarray*}
However, since the inequality
$\tilde{\varphi}_{-1}^2 \circ f^{-1}(x)-g\circ f^{-1}(x)> 0$ for all $x\in [x_{-2}, x_{-1}]$ may not be true,
we cannot obtain
$\varphi_{-2}(x) > \xi_0$.
 We also cannot obtain a weaker condition $\varphi_{-2}(x) \ge x$ for all $x\in [x_{-2}, x_{-1}]$, i.e.,
 $$
 \tilde{\varphi}_{-1}^2(x) \ge h(f(x)) + g(x),~~~~\forall x\in [x_{-1},x_0],
 $$ 
without an additional assumption on functions $f,g$ and $h$.
This prevents us from constructing $\varphi_{-3}$.
By the inequality $\xi_0 < x_{-1}$ and \eqref{gxi0}, 
i.e., $\varphi_{-1}(x) > \xi_0$ for all $x\in [x_{-1},x_0]$, we cannot deduce
$\varphi_{-1}(x) \ge x$ for all $x\in [x_{-1},x_0]$. Therefore, we cannot use $\varphi_{-1}$ and $\varphi_*$ to define $\varphi_{-2}$,
as doing in \eqref{-2}.
On the other hand,
{\it if we do not consider (\ref{gxg}),
our construction given in the proof of Theorem~\ref{thm4}
is not applicable on $[x_0, +\infty)$.}
In fact, 
without (\ref{gxg}), we have
\begin{eqnarray}\label{lex}
g(x) < x,~~~~\forall x\in J \subset [x_1, +\infty),
\end{eqnarray}
where $J$ is an interval, i.e., function $g$ lies below the diagonal as $x\in J$.
Doing as in the proof of Theorem \ref{thm4}, 
we can construct a sequence of functions $\{\varphi_i\}_{i\ge 0}$ and 
a strictly increasing sequence of points $\{x_i\}_{i\ge 0}$ alternately 
such that $\varphi_i : [x_i,x_{i+1}] \to [x_{i+1}, x_{i+2}]$ is an increasing homeomorphism for each $i\ge 0$,
but
we cannot show $x_i \to +\infty$ as $i\to +\infty$. 
In fact, assuming that 
$x_i \to x_*$, as $i \to +\infty$,
doing as in the proof of Theorem \ref{thm4} 
we get condition \eqref{x_*}, i.e., $x_* = h(\tilde{x}) + g(x_*)$,
but this may be true for some $\tilde{x}, x_*\in (x_1, +\infty)$ under condition \eqref{lex}.

As mentioned in the beginning of section 4, Theorem~\ref{thm4} is devoted to the increasing case. 
We fail to find a strictly decreasing solution of \eqref{rose} with
strictly decreasing $f$ and strictly increasing $h$ and $g$.
In fact,
choose points $x_0,x_1,x_2,x_3$ and $x_4$ such that 
$
x_4:=h(x_2)+g(x_2) < x_2 < x_0:=f ^2(x_2) < \hat{x} < x_1 := f(x_2) < x_3 := h(f(x_2)) + g(f(x_2)),
$
where $\hat{x}$ is the fixed point of $f$, and 
two decreasing homeomorphisms $\varphi_0 : [x_2,x_0] \to [x_1, x_3]$ and 
$\varphi_1 : [x_1, x_3] \to [x_4,x_2]$.
Then, as shown in \eqref{varphi2}, define
\begin{eqnarray*}
\varphi_2(x) := h \circ \varphi_0 \circ f \circ \varphi_1^{-1}(x) + g\circ \varphi_1^{-1}(x),~~~~ \forall x\in [x_4,x_2].
\end{eqnarray*}
In order to define $\varphi_2$ well, we need $f \circ \varphi_1^{-1}([x_4,x_2]) = f([x_1,x_3])\subset [x_2,x_0]$, which is equivalent to 
\begin{eqnarray}\label{x3x3}
f(x_3) \ge x_2 
\end{eqnarray}
because function $f$ is strictly decreasing and $f(x_1)= f^2(x_2)=x_0$.
Since $g$ and $h$ are strictly increasing and $\varphi_0, \varphi_1$ and $f$ are strictly decreasing, we see that
$\varphi_2$ is strictly decreasing, implying that 
$
x_5:=\varphi_2(x_4) > \varphi_2(x_2)= h\circ \varphi_0 \circ f(x_1) +g(x_1) =h(x_1) +g(x_1) =x_3.
$ 
Next, define 
\begin{eqnarray*}
\varphi_3(x) := h \circ \varphi_1 \circ f \circ \varphi_2^{-1}(x) + g\circ \varphi_2^{-1}(x),~~~~ \forall x\in [x_3,x_5].
\end{eqnarray*}
In order to define $\varphi_3$ well, we need $f \circ \varphi_2^{-1}([x_3,x_5]) = f([x_4,x_2])\subset [x_1,x_3]$, which is equivalent to 
\begin{eqnarray}\label{x4}
f(x_4) \le x_3
\end{eqnarray}
because function $f$ is strictly decreasing and $f(x_2)=x_1$.
Similarly, we see that $\varphi_3$ is strictly decreasing, implying that
$
x_6:=\varphi_3(x_5) < \varphi_3(x_3) = h\circ \varphi_1 \circ f(x_2) +g(x_2) =h(x_2) +g(x_2) =x_4.
$
In order to define
$$
\varphi_4(x) := h \circ \tilde{\varphi}_2 \circ f \circ \varphi_3^{-1}(x) + g\circ \varphi_3^{-1}(x),~~~~ \forall x\in [x_6,x_4],
$$
where 
\begin{eqnarray*}
\tilde{\varphi}_2(x) : = \begin{cases}
\varphi_0(x),~~~~  & x\in [x_2,x_0],
\\
\varphi_2(x),~~~~  & x\in [x_4, x_2),
\end{cases}
\end{eqnarray*}
we require
$
f \circ \varphi_3^{-1}([x_6,x_4]) = f([x_3,x_5]) \subset [x_4,x_0]
$
, or equivalently say, $f(x_5) \ge x_4$, i.e.,
\begin{eqnarray}\label{x5}
f(\varphi_2(x_4)) = f( h \circ \varphi_0 \circ f(x_3) + g(x_3) ) \ge x_4.
\end{eqnarray}
Even though we can choose $x_0,x_1,x_2,x_3$ and $x_4$ such that  \eqref{x3x3} and \eqref{x4},
it is difficult to decide whether \eqref{x5} is true or not.


Theorem~\ref{thm2} 
also makes some advances even if we apply it to \eqref{linear},
a special case of equation (\ref{rose}) with 
\begin{eqnarray}
h(x) := \lambda x, ~~~f(x):= x+a,~~~g(x):=\mu x.
\label{hfg}
\end{eqnarray}
Since functions given in (\ref{hfg}) satisfy {\bf (C1)} and {\bf (C2)} with constants $K =|\lambda|$ and $\alpha=1$ 
and Lip$(g)=\beta= |\mu|$,
applying Theorem \ref{thm2} to equation \eqref{linear}, 
we obtain from \eqref{contracself2} and \eqref{contracselfU1} that 
equation (\ref{linear}) has a continuous solution if
\begin{eqnarray*}
|\lambda| > \max \{2, 2\sqrt{|\mu|}\}  ~~ \text{or}
~~ 1+|\mu| < |\lambda| \le 2,
\end{eqnarray*}
which obviously is weaker than \eqref{ZWN}, a condition obtained in \cite{ZWN}. 
Besides,
Theorem~\ref{thm3} generalizes Theorem $2$ of \cite{ZYY} from the case of linear $f, g$ and $h$ to a nonlinear case.
In fact, since functions given in (\ref{hfg}) also satisfy assumptions ${\bf (C1^{'})}$-${\bf (C3^{'})}$ with $\kappa_h =\lambda, \kappa_f=\alpha=1, \kappa_g=\mu, K=|\lambda|$ and $\beta=|\mu|$,
we can also apply Theorem \ref{thm3} to equation (\ref{linear}) and rewrite
conditions \eqref{bounded} and either \eqref{contracself1} or \eqref{contracself2} as
\begin{eqnarray}
\mu \ge -\frac{1}{4} \lambda^2,
\label{bbb}
\end{eqnarray}
and either 
\begin{eqnarray}
|\lambda| > 2 ~~{\rm and}~~ |\lambda| \ge 2\sqrt{|\mu|}
\label{111}
\end{eqnarray}
or
\begin{eqnarray}
1+|\mu| < |\lambda| \le 2.
\label{222}
\end{eqnarray}
One can check that \eqref{bbb} matched with \eqref{111} is equivalent to 
\eqref{ZYY1} and that
\eqref{bbb} matched with \eqref{222} is equivalent to \eqref{ZYY2}, implying that
Theorem~\ref{thm3} gives the same conditions as Theorem $2$ of \cite{ZYY}.
Example \ref{exam2} illustrates Theorem \ref{thm3} with a nonlinear $g$.


\end{document}